\documentclass[12pt]{article}
\usepackage[final]{epsfig}
\usepackage{graphics}
\usepackage{amsmath}
\usepackage{amsfonts}
\usepackage{latexsym}
\usepackage{amssymb}
\usepackage{graphicx}
\usepackage{epstopdf}
\newtheorem{lemma}{Lemma}

\newtheorem{theorem}[lemma]{Theorem}

\begin{document}
\newcommand{\eps}{{\varepsilon}}
\newcommand{\proofend}{$\Box$\bigskip}
\newcommand{\C}{{\mathbf C}}
\newcommand{\Q}{{\mathbf Q}}
\newcommand{\R}{{\mathbf R}}
\newcommand{\Z}{{\mathbf Z}}
\newcommand{\RP}{{\mathbf {RP}}}

\title {A proof of Culter's theorem on the existence of periodic orbits in  polygonal outer billiards}
\author{Serge Tabachnikov\thanks{
Department of Mathematics,
Pennsylvania State University, University Park, PA 16802, USA;
e-mail: \tt{tabachni@math.psu.edu}
}
\\
}
\date{\today}
\maketitle
\begin{abstract}
We discuss a recent result by C. Culter: every polygonal outer billiard has  a periodic trajectory. 
\end{abstract}

This note is an exposition of a theorem proved by Chris Culter, then an undergraduate student; he obtained this result as a participant of the 2004 Penn State REU program.\footnote{The program was supported by an NSF grant; the problem solved by Culter was proposed by the author of this note.} A complete account of Culter's work  involving a more general class of maps and a more detailed analysis of their periodic orbits 
will appear in his paper, currently in progress.  

\begin{figure}[hbtp]
\centering
\includegraphics[width=2in]{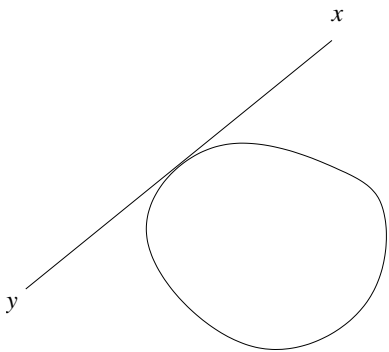}
\caption{Definition of the outer billiard map}
\label{deffig}
\end{figure}

An outer billiard table is a compact convex domain $P$. Pick a point $x$ outside $P$. There are two support lines from $x$ to $P$; choose  one of them, say, the right one from the view-point of $x$, and reflect $x$ in the support point. One obtains a new point, $y$, and the transformation $T: x\mapsto y$ is the outer (a.k.a. dual) billiard map, see figure \ref{deffig}. The map $T$ is not defined if the support line has a segment in common with the outer billiard table. In this note, $P$ is a convex $n$-gon; the set of points for which $T$ or any of its iterations is not defined is contained in a countable union of lines and has zero measure. For ease of exposition, we assume that $P$ has no parallel sides.

Outer billiards were introduced in \cite{N} and popularized in \cite{Mo1,Mo2}; we refer to \cite{D-T,Ta0,Ta6} for surveys. Here we are concerned with the existence of periodic trajectories of the outer billiard map. For the conventional, inner, billiards it is an outstanding open problem whether every polygon has a periodic billiard path. The best result so far is a theorem of R. Schwartz: every obtuse triangle with the obtuse angle not greater than $100^{\circ}$ has a periodic trajectory, see \cite{Sc1,Sc2,Sc3}. Note also that both inner and outer polygonal billiards on the sphere $S^2$ may have no periodic trajectories at all, see \cite{G-T}.

It will be convenient to consider the second iteration $T^2$ of the outer billiard map. 
Connecting the consecutive points of a periodic  trajectory of $T^2$, one obtains a closed polygonal line. The number of turns made by this line about the billiard table is called the rotation number. 
The main result is as follows.

\begin{theorem} \label{main}
The  map $T^2$ has a periodic trajectory that lies outside of any compact neighborhood of $P$ and has rotation number 1.
\end{theorem}

\paragraph {Proof of Theorem.}
For every outer billiard, not necessarily polygonal, the asymptotic dynamics of the map $T^2$ at infinity has the following description; see the sited surveys or \cite{Ta1,Ta4,Ta5}.  A bird's eye view of a outer billiard  is almost a point and the map $T$ is almost the reflection in this point. More precisely, after rescaling, the distance between a point $x$ and  $T^2(x)$ is very small, and the
evolution of a point under $T^2$ appears a continuous clockwise motion along a  centrally symmetric  curve $R$. 

\begin{figure}[hbtp]
\centering
\includegraphics[width=3in]{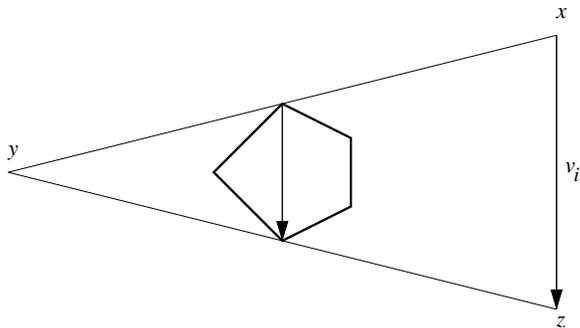}
\caption{The second iteration of the outer billiard map}
\label{twice}
\end{figure}

In our case, $R$ is a convex $2n$-gon,  and each vector $(x, T^2(x))$ belongs to a finite set  $\{\pm v_1,\dots, \pm v_n\}$. These vectors are as follows. For every direction, other than the directions of the sides, there exists a pair of parallel support lines to $P$; the vector $v_i$ is twice the vector connecting the respective support vertices of $P$, see figure \ref{twice}. For example, if $P$ is a triangle then $R$ is an affine-regular hexagon.

Consider the polygon $R$ (defined up to dilation). To every side of $R$ there corresponds ``time", the ratio of the length of this side to the magnitude of the respective vector $v_i$. One obtains a collection of ``times" $(t_1,\dots,t_k)$, defined up to a common factor. 
The polygon $P$ is called quasi-rational if all these numbers are rational multiples of each other. For example, lattice polygons are quasi-rational and so are affine-regular ones. It is known that the orbits of the outer billiard about a quasi-rational polygon are bounded, see \cite{G-S,Ko,S-V} or the cited surveys. Recently R. Schwartz proved that polygonal outer billiards may have orbits escaping to infinity  \cite{Sc4,Sc5}.

\begin{figure}[hbtp]
\centering
\includegraphics[width=3.8in]{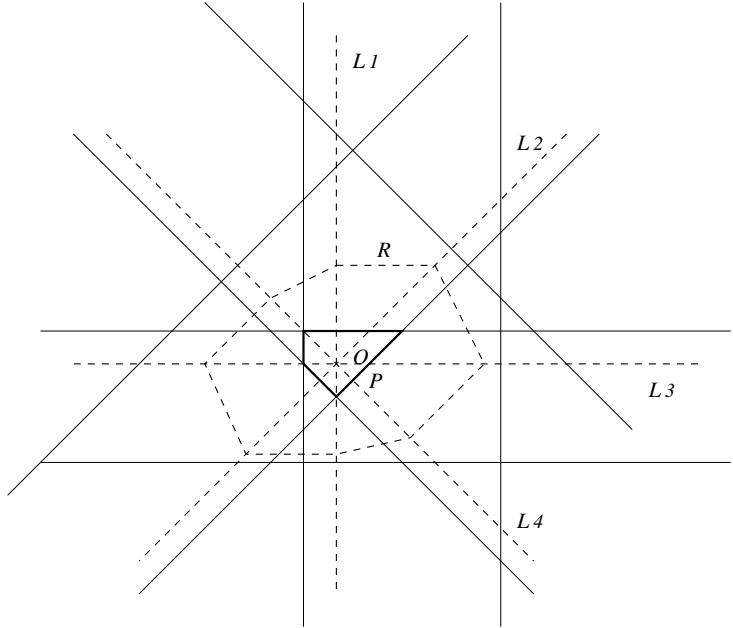}
\caption{The lines $L_i$, the strips $S_i$ and the polygon $R$}
\label{strips}
\end{figure}

The actual map $T^2$, sufficiently far away from $P$, is a piece-wise parallel translation through the vectors $\pm v_1,\dots, \pm v_n$. The discontinuities are $2n$ rays: the clockwise extensions of the sides of $P$ and the reflections of these rays in the opposite vertices of $P$ (a vertex opposite to a side is the one farthest from it). The lines containing these $2n$ rays form $n$ strips $S_1,\dots,S_n$ whose intersection contains $P$, see figure \ref{strips}.

Choose an origin $O$ inside $P$ and consider the lines $L_1,\dots, L_n$ through $O$ parallel to the sides  of $P$. Fix the above described polygon $R$ so that $O$ is its center. Denote by $qR$  the dilation of $R$ with coefficient $q$. These polygons can be constructed by choosing a starting point on $L_1$, drawing the line in the direction $v_1$ until its intersection with $L_2$, then drawing the line  in the direction $v_2$ until its intersection with $L_3$, etc. 
 
Let $p_1,\dots,p_n$ be positive integers. Denote by $Q(p_1,\dots,p_n)$ the centrally symmetric $2n$-gon whose sides are given by the vectors 
$$p_1 v_1,p_2 v_2,\dots, p_nv_n,-p_1v_1,\dots,-p_n v_n$$
 and whose center is $O$. We wish to show that, for an appropriate choice of $p_1,\dots,p_n$, the polygon $Q(p_1,\dots,p_n)$ is an orbit of the map $T^2$. For this, the vertices of $Q(p_1,\dots,p_n)$ should lie inside the strips $S_i$ (the opposite vertices in the same strip).

Clearly, there is  $\eps>0$ (depending only on $P$ and the choice of the origin) such that if the vertices of an $2n$-gon $Q$ are $\eps$-close to the respective vertices of a polygon $qR$ then the vertices of $Q$ lie inside the strips $S_i$. We claim that there exist arbitrarily large real $q$ and integers $p_1,\dots,p_n$ such that the respective vertices of  $qR$ and $Q(p_1,\dots,p_n)$ are within $\eps$ from each other. 

For the claim to hold, it will suffice to have
\begin{equation} \label{appr}
|q t_i - p_i| < \delta,\ \ i=1,\dots,n
\end{equation}
where $\delta>0$ is a small enough constant. Indeed, the first vertex of the polygon $qR$ is  
$$
 -\frac{1}{2} \sum_1^n q t_i v_i, 
 $$
whereas that of the polygon $Q(p_1,\dots,p_n)$ is
$$
 -\frac{1}{2} \sum_1^n p_i v_i,
$$
and similarly  for the other vertices. 

Finally, consider the torus $T^n=\R^n/\Z^n$, and let $F_t$ be the constant flow with the 
vector $(t_1,\dots,t_n)$. Then (\ref{appr}) means that $F_q(O)$ is $\delta$-close to $O$ where $O=(0,\dots,0)$. Indeed, the flow $F_t$ is either periodic, and then $F_q(O)=O$ for $q$ forming an arithmetic progression, or quasi-periodic and thus returning arbitrarily close to the initial point infinitely often. 
\proofend

\paragraph {\bf Remarks.}

 1. A composition of a number of central symmetries  is either a central symmetry or a parallel translation. It follows that a $k$-periodic point of the outer billiard map about a polygon has a polygonal neighborhood consisting of periodic points with period $k$ or $2k$ (the latter holds if $k$ is odd).

2. The density of the numbers $q$ satisfying (\ref{appr}) is positive. One can deduce  that the lower density of the set of periodic trajectories described in Theorem \ref{main} is also positive.

3. A periodic trajectory of the polygonal outer billiard map is called stable if, under  an arbitrary small perturbation of the outer billiard polygon $P$, the trajectory is also perturbed but not destroyed.  A criterion for stability is known, see \cite{Ta0}. Enumerate the vertices of $P$ counterclockwise as $A_1,\dots, A_n$. An even-periodic orbit of the dual billiard map is encoded by the sequence vertices in which the consecutive reflections occur. One obtains a cyclic word $W$ in the letters $A_1,\dots, A_n$. The orbit is stable if and only if each appearance of every letter in an odd position in $W$ is balanced by its appearance in an even position. By this criterion, the periodic trajectories  of Theorem \ref{main}  are stable.

\bigskip

{\bf Acknowledgments}. Many thanks to C. Culter for numerous discussions and to R. Schwartz for his interest.  The  author was partially supported by an NSF grant DMS-0555803.

\end{document}